\documentclass[11pt]{article}
\usepackage{amsmath,amssymb,url,amsthm}

\newtheorem{theorem}{Theorem}
\newtheorem{lemma}[theorem]{Lemma}

\newtheorem*{conjecture}{Conjecture}
\newtheorem*{theorem*}{Theorem}

\newcommand{\bR}{\mathbb{R}}
\newcommand{\Var}{\mathrm{Var}}
\newcommand{\Bin}{\mathrm{Bin}}

\begin{document}
\title{Connectivity of the Uniform Random Intersection Graph}
\author{Simon R. Blackburn\thanks{Corresponding author. This research was partially supported by E.P.S.R.C.\  grants EP/D053285/1 and EP/F056486/1.}\,\,~and Stefanie Gerke\\
Department of Mathematics\\
Royal Holloway, University of London\\
Egham, Surrey TW20 0EX, United Kingdom\\
{\tt s.blackburn@rhul.ac.uk, stefanie.gerke@rhul.ac.uk}}
\maketitle

\begin{abstract}
A \emph{uniform random intersection graph} $G(n,m,k)$ is a random
graph constructed as follows. Label each of $n$ nodes by a randomly
chosen set of $k$ distinct colours taken from some finite set of
possible colours of size $m$. Nodes are joined by an edge if and only
if some colour appears in both their labels. These graphs arise in the
study of the security of wireless sensor networks, in particular when
modelling the network graph of the well known key predistribution
technique due to Eschenauer and Gligor.

The paper determines the threshold for connectivity of the graph
$G(n,m,k)$ when $n\rightarrow \infty$ in many situations. For example,
when $k$ is a function of $n$ such that $k\geq 2$ and $m=\lfloor
n^\alpha\rfloor$ for some fixed positive real number $\alpha$ then
$G(n,m,k)$ is almost surely connected when
\[
\liminf k^2n/m\log n>1,
\]
and $G(n,m,k)$ is almost surely disconnected when
\[
\limsup k^2n/m\log n<1.
\]
\end{abstract}


\newpage

\section{Introduction}

\subsection{Notation and motivation}

The \emph{uniform random intersection graph} $G(n,m,k)$ is a random
graph defined as follows. Let $V$ be a set of $n$ nodes, and let $M$
be a set of $m$ colours. To each node $v\in V$ we assign a subset
$F_v\subseteq M$ of $k$ distinct colours, chosen uniformly and
independently at random from the $k$-subsets of $M$. We join distinct
nodes $u,v\in V$ by an edge if and only if $F_u\cap
F_v\not=\emptyset$. This paper studies the connectivity threshold of
uniform random intersection graphs.

The study of $G(n,m,k)$ is motivated by an application to wireless
sensor networks (WSNs). A WSN is a collection of (usually very small)
sensor devices that are able to communicate wirelessly. Sample
applications where WSNs might be used include disaster recovery,
wildlife monitoring and military situations. Sensors' computational
abilities are assumed to be severely limited by their size and battery
life. The sensor network is designed to be deployed in an unstructured
environment (sensors might be scattered from an aeroplane, for
example). On deployment the individual sensors need to form a secure
wireless network that is connected, but should also be robust against
the compromise of individual sensor's secret data due to malfunction
or capture. The classic WSN technique to accomplish this is due to
Eschenauer and Gligor~\cite{Eschenauer}: each sensor is preloaded with
$k$ distinct encryption keys, randomly taken from a pool of $m$
possible keys. Two sensors can form a secure link if they are within
wireless communication range and they share one or more encryption
keys. The uniform random intersection graph models this situation in
the case when all sensors are within communication range. (In the
terminology of the subject, a uniform random intersection graph is a
\emph{network graph} for Eschenauer--Gligor key predistribution).

The application requires the network to be connected with high
probability. Looking at other results in random graph theory, we would
expect the parameters $n$, $m$ and $k$ to exhibit a threshold
behaviour with respect to connectivity: for most parameters we would
expect that the probability that $G(m,n,k)$ is connected is either
very high or very low. It is important to understand the connectivity
threshold (the area of the parameter space bordering the regions of
low and high connectivity probability) as precisely as possible, as
this threshold effects the choice of parameters in the
Eschenauer--Gligor scheme. Eschenauer and Gligor, and most of the
subsequent WSN literature, model the uniform random intersection graph
as a classical Erd\H{o}s--R\'enyi random graph $G(n,p)$, a graph with
$n$ vertices whose edges are chosen randomly and independently with a
fixed probability $p$. They then use the asymptotic behaviour of
Erd\H{o}s--Renyi random graphs to find good parameters for the
scheme. For distinct nodes $u,v\in G(n,m,k)$, the probability that
$uv$ is an edge is $p$, where
\[
p=1-\frac{\binom{m-k}{k}}{\binom{m}{k}}\approx \frac{k^2}{m}.
\]
(To see why this approximation holds, note that $u$ is assigned $k$
colours and the probability that each colour is assigned to $v$ is
$k/m$.) So the WSN literature models $G(n,m,k)$ by the
Erd\H{o}s--R\'enyi random graph $G(n,p)$ where $p=k^2/m$. It is well
known that the connectivity threshold of $G(n,p)$ occurs when
$p\approx (\log n)/n$. So modelling $G(n,m,k)$ as an
Erd\H{o}s--R\'enyi random graph predicts that the connectivity
threshold lies at the point when $k^2/m\approx (\log n)/n$. Though
simulations support this threshold, modelling $G(n,m,k)$ in this way
is unsatisfactory since the behaviour of $G(n,p)$ and $G(n,m,k)$ is
sometimes radically different. For example, we expect many more
triangles in $G(n,m,k)$ than in $G(n,p)$, especially when $k$ is
small. (When $u,v,w\in G(n,m,2)$ are distinct vertices such that $uv$
and $vw$ are edges, then the probability that $uw$ is an edge is more
than $1/2$, since this is the probability that $v$ shares the same
colour with both $u$ and $w$.)

\subsection{Our results}

Let $k$ and $m$ be functions of $n$. Our proof techniques and results
depend heavily on whether $m\geq n$ or not, so we discuss these
two cases separately.

Suppose that $m\geq n$. We will show
(Theorem~\ref{thm:urig_connectedbigm}) that $G(n,m,k)$ is
asymptotically almost surely connected when
$\lim\inf_{n\rightarrow\infty} k^2n/(m\log n)>1$. (By an event
occurring asymptotically almost surely, we mean that the probability of
the event tends to $1$ as $n\rightarrow\infty$.) This threshold is
tight: we will show that $G(n,m,k)$ is asymptotically almost surely
disconnected when $\lim\sup_{n\rightarrow\infty} k^2n/(m\log n)<1$. Di
Pietro, Mancini, Mei, Panconesi and
Radhakrishnan~\cite{DiPietro,DiPietronew} give a weaker form of
Theorem~\ref{thm:urig_connectedbigm}: that $G(n,m,k)$ is almost surely
connected when $\lim\inf_{n\rightarrow\infty} k^2n/(m\log n)>8$. (The
journal version of their paper~\cite{DiPietronew} only claims that
$G(n,m,k)$ is almost surely connected when
$\lim\inf_{n\rightarrow\infty} k^2n/(m\log n)>17$.) Di Pietro et al
also observe that $G(n,m,k)$ is almost surely disconnected when
$k^2n/(m\log n)\rightarrow 0$ as $n\rightarrow\infty$.  Part of our
proof of Theorem~\ref{thm:urig_connectedbigm} is inspired by their
techniques. We comment that there is a gap we are unable to bridge in their
proof, which means that we take a subtly different approach to theirs:
we discuss this at the end of Section~\ref{sec:connectedbigm}.

We now turn to the case when $m\leq n$. We show (see
Section~\ref{sec:connectedk2proof}) that whenever $(4n/m)-\log
n\rightarrow \infty$ as $n\rightarrow \infty$ then $G(n,m,k)$ is
asymptotically almost surely connected. We will show (see
Theorem~\ref{thm:urig_connectedk2} below) that this threshold is tight
in the case when $k=2$. This settles the case, for example, when
$m=o(n/\log n)$. We note that this case is also a consequence of
recent work of Godehardt, Jaworski and Rybarczyk~\cite{GJR}, who show
that when $k$ is fixed, $G(n,m,k)$ is asymptotically almost surely
connected whenever $n$ is a function of $m$ such that $(kn/m)-\log
m\rightarrow\infty$ as $m\rightarrow\infty$. We believe that their
result is not tight: see Section~\ref{sec:openproblems} for a
discussion.

This leaves a narrow range of parameters not covered by our results,
when $m$ grows just a little more slowly than $n$. Though this range is too
small to be of significance in applications, there are some interesting
mathematical questions here. We comment on this in the final
section of the paper. By constraining $m$ to be of the form $m=\lfloor
n^\alpha\rfloor$ where $\alpha$ is a fixed positive real number, we
avoid this gap and obtain the following easy to state
summary of our results:

\begin{theorem}
\label{thm:urig_connectedalpha}
Let $\alpha\in\bR$ be positive.  Let $k$ and $m$ be functions of $n$
such that $k\geq 2$ and $m=\lfloor n^\alpha\rfloor$.
\begin{itemize}
\item[\emph{(i)}] Suppose that
\begin{equation}
\label{eqn:infassumptionalpha}
\liminf_{n\rightarrow\infty} \frac{k^2n}{m\log n}>1.
\end{equation}
Then asymptotically almost surely $G(n,m,k)$ is connected.
\item[\emph{(ii)}] Suppose that
\[
\limsup_{n\rightarrow\infty} \frac{k^2n}{m\log n}<1.
\]
Then asymptotically almost surely $G(n,m,k)$ is not connected.
\end{itemize}
\end{theorem}
\noindent

\subsection{Related results}

Other properties of $G(n,m,k)$ besides connectivity have been
studied. For example, Godehardt and Jaworski~\cite{GJ} have results on
the distribution of the number of isolated vertices of $G(n,m,k)$ when
$nk^2/m\log n$ tends to a constant; Bloznelis, Jaworski and
Rybarczyk~\cite{BJR} determine the emergence of the giant component
when $n(\log n)^2=o(m)$; Jaworski, Karo\'nski and Stark~\cite{JKS}
study the vertex degree distribution of random intersection graphs.

A related, non-uniform, definition of a random intersection graph has
been studied as part of the modelling of clustering in real-world
networks (see~\cite{Behrisch,Fill,SingerCohen,Stark}, for example). We
define the (non-uniform) \emph{random intersection graph} $G(n,m,p)$
exactly as in the definition of $G(n,m,k)$ above, except we choose the
subsets $F_v$ differently: each $F_v$ is constructed by the rule that
each colour $c\in M$ lies in $F_v$ independently with probability
$p$. (Thus the set $F_v$ is likely to vary in size as $v$ varies, and
will have expected size $pm$.) In her thesis,
Singer-Cohen~\cite{SingerCohen} establishes connectivity thresholds
for $G(n,m,p)$. To compare her results with
Theorem~\ref{thm:urig_connectedalpha}, consider the case when $p=k/m$,
so the expected size of a set $F_v$ is $k$. When $\alpha>1$,
Singer-Cohen shows that the connectivity threshold lies at
$p=\sqrt{(\log n) / nm}$, which agrees with the threshold of
Theorem~\ref{thm:urig_connectedalpha} (though Singer-Cohen's threshold
is sharper). In fact, when $m$ is large compared to $n$ this agreement
is a consequence of standard concentration results. When
$\alpha\leq 1$, Singer-Cohen shows that the connectivity threshold
lies at $p=\log n /m$, which is much higher than the threshold of
Theorem~\ref{thm:urig_connectedalpha}. The intuition here is that when
$m$ is small there are some nodes $v$ in $G(n,m,p)$ with $F_v$ much
smaller than $pm$ (indeed, $F_v$ may even be empty). It is these nodes
that provide the dominant obstacle to connectivity in $G(n,m,p)$ when
$\alpha\leq 1$. This also shows that $G(n,m,p)$ may behave differently
to $G(n,m,k)$.

\subsection{The structure of the paper}

The remainder of the paper is structured as
follows. Section~\ref{sec:isolatedproof} establishes the threshold for
the existence of isolated vertices in $G(n,m,k)$, using the first and
second moment methods; this result is sufficient to establish
Theorem~\ref{thm:urig_connectedalpha}~(ii). Section~\ref{sec:connectedk2proof}
specialises to the case when $k=2$, and proves
Theorem~\ref{thm:urig_connectedalpha}~(i) when $\alpha<1$.
Section~\ref{sec:connectedbigm} proves
Theorem~\ref{thm:urig_connectedalpha}~(i) when $\alpha\geq
1$. Finally, Section~\ref{sec:openproblems} discusses prospects of
establishing tighter connectivity thresholds for $G(n,m,k)$.

\section{Isolated vertices}
\label{sec:isolatedproof}

We aim to prove the following theorem on the probability of an
isolated vertex appearing in $G(n,m,k)$.

\begin{theorem}
\label{thm:urig_isolated}
Let $k$ and $m$ be functions of $n$.
\begin{itemize}
\item[\emph{(i)}] Suppose that
\begin{equation}
\label{eqn:big_parameters}
\frac{k^2n}{m}=(\log n)+\omega
\end{equation}
where $\omega\rightarrow\infty$ as $n\rightarrow\infty$. Then
almost surely $G(n,m,k)$ does not contain an isolated vertex.
\item[\emph{(ii)}] Suppose that
\begin{equation}
\label{eqn:small_parameters}
\frac{k^2n}{m}=(\log n)-\omega
\end{equation}
where $\omega\rightarrow\infty$ as $n\rightarrow\infty$. Then
almost surely $G(n,m,k)$ contains an isolated vertex.
\end{itemize}
\end{theorem}

The proof of this theorem is an application of standard techniques
from random graph theory: we include the proof for completeness. We
remark that Godehardt and Jaworski have much stronger results on the
distribution of the number of isolated vertices on the threshold: in
particular, they determine the distribution when $(k^2n/m) - \log n
\rightarrow c$ for some constant $c$; see~\cite{GJ} for a statement of
their results. Note that (in contrast to many situations in random
graph theory) it is not at all clear that Theorem~\ref{thm:urig_isolated}
immediately follows from their result: problems occur with a reduction
as, for example, $k$ has to be integer and if one changes $k$ by $1$
then $k^2n/m$ may vary by a factor greater than $\log n$.
\begin{proof}
For $v\in V$, define the random variable $X_v$ by
\[
X_v=\left\{\begin{array}{cl}
1&\text{ if $v$ is isolated,}\\
0&\text{ otherwise.}
\end{array}\right.
\]
Define $X=\sum_{v\in V} X_v$. So $E(X)$ is the expected number of
isolated vertices in $G(n,m,k)$. Note that, by linearity of
expectation, $E(X)=nE(X_u)$, where $u\in V$ is any fixed
vertex. A vertex is isolated if and only if $F_v\cap
F_u=\emptyset$ for all $v\in V\setminus\{u\}$. Hence
\begin{align*}
E(X)&=n\left(
\frac{\binom{m-k}{k}}{\binom{m}{k}}\right)^{n-1}=n\left(\prod_{i=0}^{k-1}\frac{m-k-i}{m-i}\right)^{n-1}\\
&=n\left(\prod_{i=0}^{k-1}
1-\frac{k}{m-i}\right)^{n-1}.
\end{align*}

Suppose that~\eqref{eqn:big_parameters} holds. We show that then
$E(X)\rightarrow 0$ and the result follows by Markov's inequality.  We
have that
\begin{align*}
E(X)&\leq n\left(1-\frac{k}{m}\right)^{k(n-1)} \leq
n\exp\left(-\frac{k^2(n-1)}{m}\right)\\
&=\exp\left( -w+ o(w) \right)
\text{ by~\eqref{eqn:big_parameters}.}
\end{align*}
So $E(X)\rightarrow 0$, as required.

We now aim to prove Part~(ii) of the theorem using the second moment method.  Note first that \eqref{eqn:small_parameters} implies that $k=o(m)$, and thus for sufficiently large $n$
\[ \frac{k}{m-k}\sqrt{k(n-1)}\leq \frac{2 k^2n}{m\sqrt{n}} =o(1).\]
Since $(1-p)^x=\exp({-px+o(1)})$ whenever $p\sqrt{x}=o(1)$ we have
\begin{align*}
E(X)&= n\left(\prod_{i=0}^{k-1} 1-\frac{k}{m-i}\right)^{n-1}
\geq n\left(1-\frac{k}{m-k}\right)^{k(n-1)}\\
&= n\exp\left(-\frac{k^2n}{m-k}+o(1)\right)= n\exp\left(-\frac{k^2n}{m}+o(1)\right)\\
&= \exp(w+o(1))
\end{align*}
which tends to infinity as $n\rightarrow\infty$. The second moment method  now implies the result we require, provided that
we can show that $\Var(X)\ll E(X)^2$. Now
\[
\Var(X)=E(X^2)-E(X)^2\geq 0,
\]
and so it suffices to show that $E(X^2)=(1+o(1))E(X)^2$. Note that
\[
E(X^2)=E(X)+n(n-1)E(X_{u_1}X_{u_2}),
\]
where $u_1,u_2$ are fixed vertices. Since $E(X)\rightarrow\infty$, it therefore
suffices to prove that
\begin{equation}
\label{eqn:varaim}
\frac{n(n-1)E(X_{u_1}X_{x_2})}{E(X)^2}\rightarrow 1 \text{ as
}n\rightarrow\infty.
\end{equation}
Note that $X_{u_1}X_{u_2}$ takes the value $1$ exactly when $u_1$ and
$u_2$ are both isolated. For $u_1$ and $u_2$ to both be isolated,
$F_{u_1}$ and $F_{u_2}$ should be disjoint (so there is no edge
between $u_1$ and $u_2$) and for all $v\in V\setminus\{u_1,u_2\}$ we
must have that $F_v$ is disjoint from $F_{u_1}\cup F_{u_2}$ (so there
is no edge from $v$ to either of $u_1$ or $u_2$). Thus
\begin{align*}
E(X_{u_1}X_{u_2}) &=
\frac{\binom{m-k}{k}}{\binom{m}{k}}\left(\frac{\binom{m-2k}{k}}{\binom{m}{k}}\right)^{n-2}\\
&=\frac{\binom{m-k}{k}}{\binom{m}{k}}\left(\frac{\binom{m-2k}{k}}{\binom{m}{k}}\right)^{-2}\left(\frac{\binom{m-2k}{k}}{\binom{m}{k}}\right)^{n}\\
&=\exp\left(-\frac{2k^2n}{m}+o(1)\right)
\end{align*}
as before. Since we proved above that
\[
E(X)=n\exp\left(-\frac{k^2n}{m}+o(1)\right),
\]
we see that~\eqref{eqn:varaim} holds, as required.
\end{proof}

\section{The case when $k=2$ or $m=o(n/\log n)$}
\label{sec:connectedk2proof}

In this section we prove the following theorem concerning the case
when each vertex is assigned a set of colours of size two.
\begin{theorem}
\label{thm:urig_connectedk2}
Let $m$ be a function of $n$.
\begin{itemize}
\item[\emph{(i)}] Suppose that
\begin{equation}
\label{eqn:big_parametersk2}
\frac{4n}{m}=(\log n)+\omega
\end{equation}
where $\omega\rightarrow\infty$ as $n\rightarrow\infty$. Then
almost surely $G(n,m,2)$ is connected.
\item[\emph{(ii)}] Suppose that
\begin{equation*}
\frac{4n}{m}=(\log n)-\omega
\end{equation*}
where $\omega\rightarrow\infty$ as $n\rightarrow\infty$. Then
almost surely $G(n,m,2)$ is not connected.
\end{itemize}
\end{theorem}
We remark that this theorem implies that $G(n,m,k)$ is asymptotically
almost surely connected whenever $m=o(n/\log n)$ (and, in particular,
Theorem~\ref{thm:urig_connectedk2} implies
Theorem~\ref{thm:urig_connectedalpha} holds when $\alpha<1$). To see
this, we first choose $2$ colours for each vertex from the $m$
available colours uniformly at random to obtain an instance of
$G(n,m,2)$. As $m=o(4n/\log n)$ we have $\log n=o(4n/m)$ and thus by
Theorem~\ref{thm:urig_connectedk2} the graph $G(n,m,2)$ is
asymptotically almost surely connected. If we now choose $k-2$ more
colours for each vertex from the remaining available colours uniformly
at random then each vertex has been assigned $k$ colours uniformly at
random, and so we have obtained an instance of $G(n,m,k)$. Moreover
the newly chosen colours can only add edges to the graph and thus the
instance of $G(n,m,k)$ is more likely to be connected than the
instance of $G(n,m,2)$.

To prove Theorem~\ref{thm:urig_connectedk2} we first prove the
following lemma which says that we only have to consider values of $m$
that are not too small compared with~$n$.

\begin{lemma}
\label{lem:restrictionk2}
It is sufficient to prove Part~\emph{(i)} of
Theorem~\emph{\ref{thm:urig_connectedk2}} in the case when
\begin{equation}
\label{eqn:restriction_assumptionk2}
\frac{n}{m\log n}\leq 1.
\end{equation}
\end{lemma}
\begin{proof} Suppose that we have proved Part~(i) of
Theorem~\ref{thm:urig_connectedk2} under the additional
assumption~\eqref{eqn:restriction_assumptionk2}. Suppose
that~\eqref{eqn:restriction_assumptionk2} is not satisfied. To prove the
lemma, it is sufficient to show that we may replace $m$ by a larger
function $m'$ of $n$ such that $ \frac{4n}{m'}-\log n \rightarrow  \infty$ and $\frac{4n}{m'\log n}\leq 4$, with the
property that $G(n,m',2)$ is less likely to be connected than
$G(n,m,2)$.

Define $m'$ by setting $m'=m$
whenever~\eqref{eqn:restriction_assumptionk2} is satisfied; otherwise
let $\ell$ be the unique positive integer such that
\[
2\leq \frac{4n}{2^\ell m\log n}\leq 4
\]
and define $m'=2^\ell m$.  Note that
\[
\frac{4n}{m'}-\log n\rightarrow\infty
\]
as $n\rightarrow\infty$ since whenever
$m\not=m'$ we have that
\[
\frac{4n}{m'}=\frac{4n}{2^\ell m\log n}\log n\geq 2\log n,
\]
by our choice of $\ell$.

It remains to show that $G(n,m',2)$ is less likely to be connected
than $G(n,m,2)$.

Let $M'$ be a set of $m'$ colours. Partition $M'$ into $m$ classes,
each of size $2^\ell$. Identify the set $M$ of $m$ colours with the
classes of this partition. We generate an instance of $G(n,m,2)$ as
follows. Firstly, we generate an instance of $G(n,m',2)$, so each node
$v$ is assigned a set $F'_v\subseteq M'$ of size~$2$. Secondly, by
replacing each colour by the class containing it we assign a set of at
most $2$ colours from $M$ to each vertex. Thirdly, for those vertices
assigned only one  colour from $M$,  we assign an additional
colour uniformly and independently at random. Note that this process
does indeed generate an instance of $G(n,m,2)$, since the vertices
assigned one colour from $M$ in the second step are coloured
uniformly and independently. To
see that $G(n,m,2)$ is more likely to be connected than $G(n,m',2)$,
note that each of the last two steps adds edges to the graph (where
the adjacency relation of the graph at the end of the second step is
chosen to be the obvious one).
\end{proof}

\noindent\emph{Proof of Theorem~\emph{\ref{thm:urig_connectedk2}}.}
Part~(ii) of Theorem~\ref{thm:urig_connectedk2} follows from Part~(ii) of
Theorem~\ref{thm:urig_isolated}, since a graph with an isolated vertex
cannot be connected. So it suffices to prove Part~(i) of the theorem. Moreover by
Lemma~\ref{lem:restrictionk2} we may assume for the remainder of the proof that
\begin{equation}\label{eq_lowerk2}
\frac{4n}{m\log n}\leq 4.
\end{equation}

Given a graph $G(n,m,2)$, we define the corresponding \emph{colour
graph} $H(n,m,2)$ as follows. The vertices of $H(n,m,2)$ are the set
$M$ of colours. Two distinct vertices $x$ and $y$ of $H(n,m,2)$ are
connected by an edge if and only if some vertex $v$ in $G(n,m,2)$ is
assigned the set $\{x,y\}$ of colours (in other words, if there exists
$v\in G(n,m,2)$ such that $F_v=\{x,y\}$). Thus $H(n,m,2)$ has $m$
vertices and at most $n$ edges.

We claim  that the colour graph $H(n,m,2)$ asymptotically almost surely contains at least $n-(\log
n)^3$ edges. To prove the claim we define for any two distinct vertices $u,v\in G(n,m,k)$, a random variable $X_{u,v}$ by
\[
X_{u,v}=\left\{\begin{array}{cl}
0&\text{ if }F_u\not=F_v,\\
1&\text{ if }F_u=F_v,
\end{array}\right.
\]
and let $X=\sum X_{u,v}$, where the sum is over all pairs of distinct
vertices in $G(n,m,2)$. Now $E(X_{u,v})=\binom{m}{2}^{-1}$, and
so~\eqref{eq_lowerk2} and linearity of expectation imply that
\[  E(X)=\binom{n}{2}\binom{m}{2}^{-1} \leq \frac{2n^2}{m^2} \leq 2(\log n)^2. \]
Markov's inequality now implies that
\[
\Pr\big(X\geq (\log n)^3\big)\leq 2(\log n)^2/(\log n)^3=2(\log n)^{-1}\rightarrow 0,
\]
and so asymptotically almost surely there are at most $(\log n)^3$ pairs $u,v$ of
vertices such that $F_u=F_v$. When $H(n,m,2)$ has $n-i$ edges, there
must be at least $i$ pairs $u,v\in G(n,m,2)$ with $F_u=F_v$. So the
claim follows.

We say a graph is \emph{near connected} if it consists of a connected
component together with a (possibly empty) set of isolated
vertices. Note that $G(n,m,2)$ is connected if and only if the
corresponding colour graph $H(n,m,2)$ is near connected.
We may regard the edges of $H(n,m,2)$ as being obtained by sampling
$n$ times with replacement from the set of edges of the complete graph
on $m$ vertices (with the uniform distribution).  By
 our claim asymptotically almost surely $H(n,m,2)$ contains at
least $n-(\log n)^3$ edges, and so we stop the process after we have
sampled this number of distinct edges to obtain a subgraph $H'$ of
$H(n,m,2)$. Note that $H'$ is chosen uniformly from the set of all
graphs on $m$ vertices with $n-\lfloor(\log n)^3\rfloor$ edges. Since
the property of being near connected is monotone,
Theorem~\ref{thm:urig_connectedk2} will follow if we can show that $H'$
is almost surely near connected. A random graph with $m$ vertices and
$x$ edges is near connected whenever
\begin{equation}
\label{eqn:Bollobasresult}
x\geq
\frac{m}{4}(\log m+\log\log m +\omega'),
\end{equation}
where $\omega'\rightarrow\infty$ as $m\rightarrow\infty$ (see
Bollob\'as~\cite[Page~164]{Bollobas}). Now,
\[
\log m\leq \log 4+\log n-\log\log n
\]
since $4n/m\geq \log n$ by~\eqref{eqn:big_parametersk2}. Since $m\leq
n$ whenever $m$ is sufficiently large, we find that
\[
\log n\geq \log m+\log\log n-\log 4\geq \log m+\log\log m-\log 4.
\]
If we set $x=n-\lfloor(\log n)^3\rfloor$ we see that
\begin{align*}
\frac{4x}{m}&\geq \frac{4n}{m}-\frac{4(\log n)^3}{m}\\
&=\log n +\omega - o(1) \text{ by~\eqref{eqn:big_parametersk2} and \eqref{eq_lowerk2}}\\
&\geq \log m + \log\log m + \omega + O(1).
\end{align*}
Thus $\frac{4x}{m}\geq \log m+\log\log m+\omega'$ where
$\omega'\rightarrow\infty$ as $m\rightarrow\infty$,
and therefore~\eqref{eqn:Bollobasresult} holds. So $H'$ is near connected, and
the theorem follows.\hfill $\Box$

\section{The case when $m\geq n$}
\label{sec:connectedbigm}

\begin{theorem}
\label{thm:urig_connectedbigm}
Let $k$ and $m$ be functions of $n$ such that $m\geq n$.
\begin{itemize}
\item[\emph{(i)}] Suppose that
\begin{equation}
\label{eqn:infassumptionbigm}
\liminf_{n\rightarrow\infty} \frac{k^2n}{m\log n}>1.
\end{equation}
Then asymptotically almost surely $G(n,m,k)$ is connected.
\item[\emph{(ii)}] Suppose that
\[
\limsup_{n\rightarrow\infty} \frac{k^2n}{m\log n}<1.
\]
Then asymptotically almost surely $G(n,m,k)$ is not connected.
\end{itemize}
\end{theorem}

Note that this theorem implies Theorem~\ref{thm:urig_connectedalpha}
holds in the case when $\alpha\geq 1$, and so our proof of
Theorem~\ref{thm:urig_connectedalpha} is complete once we have proved
this theorem.  As before, Part~(ii) of
Theorem~\ref{thm:urig_connectedbigm} follows from Part~(ii) of
Theorem~\ref{thm:urig_isolated}, since a graph with an isolated vertex
cannot be connected. So it suffices to prove Part~(i) of the
theorem. Our proof of Part~(i) parallels and tightens the work of Di
Pietro \emph{et al}~\cite{DiPietro}.

If $G(n,m,k)$ is not connected, it has a component $S$ of size at most
$n/2$. Lemmas~\ref{lem:nosmallcomponents},
\ref{lem:nomediumcomponents} and~\ref{lem:nolargecomponents} together
show that the probability of such a component $S$ existing tends to
$0$ as $n\rightarrow\infty$, and so the theorem will follow from these
three lemmas.

Note that \eqref{eqn:infassumptionbigm} and the fact that $m\geq n$
together imply that $k\geq \sqrt{\log n}$ for all sufficiently large
$n$. In particular, $k\rightarrow\infty$ as $n\rightarrow\infty$.

\begin{lemma}
\label{lem:nosmallcomponents}
Under the conditions of Part~\emph{(i)} of
Theorem~\emph{\ref{thm:urig_connectedbigm}}, $G(n,m,k)$ asymptotically almost
surely contains no components of size $s$, with $s\leq en^{8/9}$.
\end{lemma}
\begin{proof}
We claim that it suffices to prove the lemma under the additional
assumption that
\begin{equation}
\label{eq_kbounded}
k^2\leq \frac{4m \log n }{ n }.
\end{equation}
For suppose we have proved the lemma under this additional
assumption. Given any $k$ satisfying~\eqref{eqn:infassumptionbigm},
define $k'$ by
\[
k'=\left\{\begin{array}{cl}
k&\text{ if }k^2\leq (4m\log n)/n,\\
\lfloor \sqrt{(4m\log n)/n}\rfloor&\text{ otherwise}.
\end{array}\right.
\]
Since $2\leq k'\leq k$, we may construct an instance of $G(n,m,k)$ by first
assigning $k'$ colours to each vertex to obtain an instance of
$G(n,m,k')$, and then assigning an additional $k-k'$ colours to each
vertex to obtain an instance of $G(n,m,k)$. Assigning the additional
$k-k'$ colours can only add edges to the graph, so the probability
that $G(n,m,k)$ has no component of order at most $e n^{8/9}$ is
bounded below by the corresponding probability for $G(n,m,k')$. Since
$\lim\inf (k')^2n/m\log n>1$, the probability that $G(n,m,k')$
has no component of order at most $e n^{8/9}$ tends to $1$, by the
lemma under the additional assumption~\eqref{eq_kbounded}. So our
claim follows.

For a set $S$ of vertices of size $s$, let $A_S$ be the
event that $S$ is a component of $G(n,m,k)$. Choose a constant $0<\varepsilon<1$  such that
\begin{equation}\label{eq_eps}
(1-2\varepsilon)\frac{k^2n}{m\log n}>1
 \end{equation}
 for all sufficiently large $n$. Such a constant exists
by \eqref{eqn:infassumptionbigm}.  Let $B_S$ be the event that
fewer than $(1-\varepsilon)ks$ colours are assigned to $S$. Note that
\begin{align*}
\Pr(A_S)&=\Pr(B_S)\Pr(A_S\mid B_S)+\Pr(\overline{B_S})\Pr(A_S\mid
\overline{B_S})\\
&\leq \Pr(B_S)+\Pr(A_S\mid\overline{B_S}).
\end{align*}
First, we shall give an upper bound on $\Pr(B_S)$. There are $\binom{m}{\lfloor (1-\varepsilon)ks\rfloor}$ choices for a set
of $\lfloor (1-\varepsilon)ks\rfloor$ colours; each of the $s$ vertices
in $S$ is assigned a subset of these colours with probability
$\binom{\lfloor(1-\varepsilon)ks\rfloor}{k}/\binom{m}{k}$. So
\begin{align*}
\Pr(B_S)& \leq \binom{m}{\lfloor(1-\varepsilon)ks\rfloor}\left(\frac{\binom{\lfloor(1-\varepsilon)ks\rfloor}{k}}{\binom{m}{k}}\right)^s\\
&\leq  \left( \frac{em}{(1-\varepsilon)ks}  \right)^{(1-\varepsilon)ks}
\left(\frac{(1-\varepsilon)ks}{m}\right)^{ks}\\
&\leq e^{ks} \left(\frac{ks}{m}\right)^{\varepsilon ks}.
\end{align*}
By \eqref{eq_kbounded} and since $s\leq n^{8/9}$ and $m\geq n$, we have
\[
\frac{ks}{m}\leq \sqrt{\frac{4m\log
n}{n}}\frac{n^{8/9}}{m}
\leq 2 n^{-\frac{1}{9}}\sqrt{\log n}.
\]
Since $k\rightarrow \infty$ as $n\rightarrow \infty$ we have
$\varepsilon k \rightarrow \infty $ and thus for sufficiently large
$n$
\[
\Pr(B_S) \leq \left[
\left(\frac{e^{1/\epsilon}ks}{m}\right)^{\varepsilon k}\right]^s \leq n^{-2s }.
\]

If $B_S$ does not occur, we may find a subset $K$ of colours of size
$\lceil (1-\varepsilon)ks\rceil $ that have been assigned to $S$. For $S$ to be a
component, each of the $n-s$ vertices not in $S$ must be assigned
colours that are disjoint from $K$, and so
\begin{align*}
\Pr(A_S\mid \overline{B_S})&\leq \left(\frac{\binom{\lfloor
m-(1-\varepsilon)ks\rfloor}{k}}{\binom{m}{k}}\right)^{n-s}\leq
\left(\frac{ m-(1-\varepsilon)ks}{m}\right)^{k(n-s)}\\
&\leq \exp\left(-(1-\varepsilon)\frac{s(n-s)}{n}\frac{k^2n}{m} \right)
\leq n^{-s \frac{1-\varepsilon}{1-2\varepsilon} \frac{n-s}{n}}
\text{ by \eqref{eq_eps}}\\
&\leq n^{-(1+\varepsilon)s}
\end{align*}
for sufficiently large $n$.

The event that $G(n,m,k)$ has a component of size at most $en^{8/9}$
is bounded above by the following expression, where we sum over all
subsets $S$ of vertices of size at most $en^{8/9}$:
\begin{align*}
\sum_S \Pr(A_S)&\leq\sum_S\big(\Pr(B_S)+\Pr(A_S\mid\overline{B_S})\big)\\
&\leq
\sum_{s=1}^{en^{8/9}}\binom{n}{s}\left(n^{-2s}+n^{-(1+\varepsilon)s}\right)\\
&\leq \sum_{s=1}^{\infty}n^s 2n^{-(1+\varepsilon)s}= \frac{2}{n^{\varepsilon}-1},
\end{align*}
which tends to zero as $n$ tends to infinity.
\end{proof}

\begin{lemma}
\label{lem:nomediumcomponents}
Under the conditions of Part~\emph{(i)} of
Theorem~\emph{\ref{thm:urig_connectedbigm}}, G(n,m,k) asymptotically
almost surely contains no components of size $s$, where $en^{8/9}\leq
s\leq \min\{m/k,n/2\}$.
\end{lemma}
\begin{proof}
Just as in the proof of Lemma~\ref{lem:nosmallcomponents}, we may
assume in addition that the inequality~\eqref{eq_kbounded} holds.

For a subset $S$ of vertices of size $s$, define $C_S$ to be the
event that $S$ is assigned at most $\frac{1}{4}ks$ colours. We proceed
as in Lemma~\ref{lem:nosmallcomponents}, with the event $C_S$
replacing the event $B_S$.  So the probability that $G(n,m,k)$
contains a component of the size we are interested in is bounded above
by
\[
\sum_{S} \Pr(C_S)+\Pr(A_S\mid \overline{C_S}),
\]
where we are summing over all subsets of vertices of size $s$, where
$en^{8/9}\leq s\leq \min\{m/k,n/2\}$. We wish to prove that this sum
tends to $0$ as $n\rightarrow\infty$.

We begin by showing that
\[
\sum_{S}\Pr(C_S)\rightarrow 0
\]
as $n\rightarrow\infty$. A similar argument to that in the proof of
Lemma~\ref{lem:nosmallcomponents} shows that
\[
\sum_{S}\Pr(C_S)\leq \sum_{s=\lceil
en^{8/9}\rceil}^{\lfloor \min\{m/k,n/2\}\rfloor} \binom{n}{s}\binom{m}{\lfloor ks/4\rfloor}\left(\frac{ks}{4m}\right)^{ks}.
\]
But then
\begin{align*}
\sum_{S}\Pr(C_S)&\leq \sum_s \binom{m}{s}\binom{m}{\lfloor
ks/4\rfloor}\left(\frac{ks}{4m}\right)^{ks}\\
&\leq\sum_s
\binom{m}{\lfloor
ks/4\rfloor}^2\left(\frac{ks}{4m}\right)^{ks}\\
&\leq \sum_{s=1}^{\lfloor\min\{ m/k,n/2\}\rfloor}
\left(\frac{eks}{4m}\right)^{ks/2}.
\end{align*}
We may write the summand in this last expression in the form
$(x^x)^t$, where $x=eks/4m$ and $t=2m/e$. Since $x^x$ has no internal
maxima (just a single minimum at $x=e^{-1}$), our summand is
maximized at the extremes of its range. So our summand is bounded
above by $\mu$ where
\[
\mu=\max\left\{\left(\frac{ek}{4m}\right)^{k/2},\left(\frac{e}{4}\right)^{m/2},\left(\frac{e}{4}\right)^{nk/4}\right\}=o(n^{-1}),
\]
by~\eqref{eq_kbounded} and since $k\rightarrow \infty$. Thus
\begin{align*}
\sum_{S}\Pr(C_S)\leq((n/2)+1)\mu=o(1),
\end{align*}
as required.

The event $A_S$ requires that the colours assigned to the $n-s$
elements of $V\setminus S$ are disjoint from the colours assigned to
$S$ (for otherwise there would be edges between $V\setminus S$ and
$S$), and so if $C_S$ does not occur we see that
\begin{align*}
\Pr(A_S\mid \overline{C_S})&\leq \left(\frac{\binom{m-(ks/4)}{k}}{\binom{m}{k}}\right)^{n-s}\\
&\leq\left(1-\frac{ks}{4m}\right)^{k(n-s)}.
\end{align*}
Hence
\begin{align*}
\sum_S\Pr(A_S\mid \overline{C_S})&\leq\sum_{s=\lceil en^{8/9}\rceil}^{\lfloor \min\{m/k,n/2\}\rfloor} \binom{n}{s}\left(1-\frac{ks}{4m}\right)^{k(n-s)}\\
&\leq\sum_s\left(\frac{ne}{s}\right)^s\exp\left(-\frac{sk^2}{4m}(n-s)\right)\\
&\leq \sum_s n^{\frac{1}{9}s} \exp\left(-\frac{sk^2n}{8m}\right)\\
&\leq \sum_s n^{\frac{1}{9}s} n^{-\frac{1}{8}s}\text{ (by~\eqref{eqn:infassumptionbigm})}\\
&\leq \sum_{s=1}^\infty n^{-\frac{1}{72}s} = \frac{1}{n^\frac{1}{72}-1}
\end{align*}
which tends to $0$ as $n$ tends to infinity. \end{proof}

\begin{lemma}
\label{lem:nolargecomponents}
Under the conditions of Part~\emph{(i)} of
Theorem~\emph{\ref{thm:urig_connectedbigm}}, $G(n,m,k)$ asymptotically
almost surely contains no components of size $s$, where $m/k<s\leq
n/2$.
\end{lemma}
\begin{proof}
We need to show that the probability that $G(n,m,k)$ has
a component of size $s>m/k$, where $s\leq n/2$, tends to $0$. If
$m/k\geq n/2$ this probability is $0$, so we assume that $m/k\leq
n/2$.

Let $T$ be a set of vertices of size $\lceil m/k\rceil$. Let $D_T$ be
the event that there are at least $n/2$ vertices in $V\setminus T$
having no edges to $T$. Note that if $G(n,m,k)$ contains a component
$S$ of size $s$ where $m/k<s\leq n/2$, all the events $D_T$ where
$T\subseteq S$ occur (since $V\setminus S$ has size at least $n/2$, and
the vertices in $V\setminus S$ have no edges to $S$ and so in
particular have no edges to $T$). So the probability that $G(n,m,k)$
contains a component of size $s$ where $m/k<s\leq n/2$ is bounded
above by $\sum_{T}\Pr(D_T)$, where the sum is over all subsets
$T\subseteq V$ with $|T|=\lceil m/k\rceil$.  Let $C_T$ be the event
that $T$ is assigned $m/4$ colours or fewer. We have that
\begin{align*}
\Pr(D_T)&=\Pr(C_T)\Pr(D_T\mid C_T)+\Pr(\overline{C_T})\Pr(D_T\mid \overline{C_T})\\
&\leq \Pr(C_T)+\Pr(D_T\mid \overline{C_T}),
\end{align*}
and so it suffices to show that $\sum_T \Pr(C_T)$ and
$\sum_T\Pr(D_T\mid \overline{C_T})$ both tend to $0$ as
$n\rightarrow\infty$. Now,
\begin{align*}
\Pr(C_T)&\leq \binom{m}{\lfloor m/4\rfloor}\left(\frac{\binom{\lfloor
m/4\rfloor}{k}}{\binom{m}{k}}\right)^{\lceil m/k\rceil}\\
&\leq \left( \frac{me}{m/4}\right)^{m/4}\left(\frac{m/4}{m}\right)^{m}\\
&=(4e)^{m/4}4^{-m}.
\end{align*}
As $k\rightarrow \infty$, we may assume that $k\geq 4$ when $n$ is
sufficiently large. So
\begin{align*}
\sum_T\Pr(C_T)&\leq \binom{n}{\lceil m/k\rceil}\Pr(C_T)\\
&\leq \binom{m}{\lfloor m/4\rfloor}(4e)^{m/4}4^{-m}\\
&\leq (4e)^{m/2}4^{-m}.
\end{align*}
Since $\sqrt{4e}<4$, we see that this sum tends to $0$ as
$n\rightarrow\infty$.

Let $T$ be fixed, and let $v\in V\setminus T$. Let $E_v$ be the event
that there are no edges from $v$ to $T$. Then $\Pr(E_v)$ is equal to
the probability that the colours assigned to $v$ are disjoint from the
colours assigned to $T$. Thus
\begin{align*}
\Pr(E_v\mid \overline{C_T})\leq \frac{\binom{m-\lfloor
m/4\rfloor}{k}}{\binom{m}{k}}&\leq \left(\frac{\lceil
3m/4\rceil}{m}\right)^k\\ &\leq (4/5)^k\leq (4/5)^{\sqrt{\log n}}.
\end{align*}
Note that the events $E_v$ are independent. The event $D_T$ occurs
exactly when $n/2$ or more of the events $E_v$ occur. So, writing
$p=(4/5)^{\sqrt{\log n}}$, we find
\begin{align*}
\Pr(D_T\mid \overline{C_T})&\leq \Pr\big(\Bin\big(n-\lceil m/k\rceil,\Pr(E_v\mid
C_T)\big)\geq n/2\big)\\
&\leq \Pr\big(\Bin(n,p)\geq n/2\big)\\
&\leq \exp\left(n(\tfrac{1}{2}\log 2p +\tfrac{1}{2}\log(2(1-p)))\right)\\
&\quad\text{by
the Chernoff bound (see Bollob\'as~\cite[Page~11]{Bollobas})}\\
&\leq \exp\left(-\tfrac{1}{2}n\sqrt{\log n}\log(5/4)+O(n)\right).
\end{align*}
Thus
\begin{align*}
\sum_T\Pr(D_T\mid \overline{C_T})&\leq 2^n
\exp\left(-\tfrac{1}{2}n\sqrt{\log n}\log(5/4)+O(n)\right)\\
&=\exp\left(-\tfrac{1}{2}n\sqrt{\log
n}\log(5/4)+O(n)\right)\rightarrow 0
\end{align*}
as $n\rightarrow\infty$. So the lemma follows.
 \end{proof}

We comment that our approach subtly differs from Di Pietro \emph{et
al}~\cite{DiPietro,DiPietronew}, in the following way. Let $B$ be the
event that there exists a set $S$ of the vertices of $G(n,m,k)$ with
$|S|\leq \min\{m/k,n/2\}$ which is assigned $|S|k/4$ or fewer distinct
colours. Di~Pietro \emph{et al} show that this event occurs with
negligible probability, and then perform the rest of their analysis on
the random graph obtained from $G(n,m,k)$ under the assumption that
$B$ does not occur. The colours assigned to different vertices given
that $B$ does not occur are no longer independent, but Di~Pietro
\emph{et al} seem to assume independence in their estimates. Our
approach avoids this problem by considering the individual events
$B_S$ for a fixed subset $S$ of vertices (see the proof of
Lemma~\ref{lem:nosmallcomponents} for example). The event $B_S$ only
depends on the colours assigned to vertices in $S$, so colours
assigned to vertices not in $S$ are still chosen independently when we
assume that $B_S$ does not occur.

\section{Discussion}
\label{sec:openproblems}

We conjecture that it is possible to prove a sharper threshold for
uniform random intersection graphs. Indeed, we believe that the
following conjecture is true.
\begin{conjecture}
\label{conj:urig_main}
Let $k$ and $m$ be functions of $n$.
\begin{itemize}
\item[\emph{(i)}] Suppose that
\[
\frac{k^2n}{m}=(\log n)+\omega
\]
where $\omega\rightarrow\infty$ as $n\rightarrow\infty$. Then
almost surely $G(n,m,k)$ is connected.
\item[\emph{(ii)}] Suppose that
\[
\frac{k^2n}{m}=(\log n)-\omega
\]
where $\omega\rightarrow\infty$ as $n\rightarrow\infty$. Then
almost surely $G(n,m,k)$ is not connected.
\end{itemize}
\end{conjecture}
\noindent
The results in this paper show that Part~(ii) of the conjecture holds
(see Theorem~\ref{thm:urig_isolated} in
Section~\ref{sec:isolatedproof} above). Moreover the full conjecture
holds in the special case when $k=2$ (by
Theorem~\ref{thm:urig_connectedk2} in
Section~\ref{sec:connectedk2proof}). To prove the full conjecture, a
natural approach would be to determine the correct generalisation to
hypergraphs of the threshold~\eqref{eqn:Bollobasresult} for the near
connectivity of graphs. This might allow a proof along the lines of
Section~\ref{sec:connectedk2proof}. However, as far as the authors are
aware, no sufficiently strong results for hypergraphs are currently
known: it would be interesting to see whether such results could be
established.

Let $p_\mathrm{conn}(n,m,k)$ be the probability that $G(n,m,k)$ is
connected.  It is easy to show that the function
$p_\mathrm{conn}(n,m,k)$ is non-decreasing in $k$. We proved a special
case of this fact in our comments below the statement of
Theorem~\ref{thm:urig_connectedk2}, and essentially the same proof
works in general. It seems reasonable to believe that
$p_\mathrm{conn}(n,m,k)$ is a non-increasing function of $m$ (so the
probability that $G(n,m+1,k)$ is connected is no larger than the
probability that $G(n,m,k)$ is connected) but we are not able to find
a proof of this. Can a proof be found?

\end{document}